\documentclass{article}
\usepackage{graphicx}
\usepackage{amsmath}
\usepackage{amsfonts}
\usepackage{amssymb}
\newtheorem{theorem}{Theorem}

\newtheorem{proposition}[theorem]{Proposition}
\newtheorem{remark}[theorem]{Remark}

\newenvironment{proof}[1][Proof]{\textbf{#1.} }{\ \rule{0.5em}{0.5em}}

\begin{document}

\date{}
\title{On Polyharmonic Interpolation}
\author{Werner Hau\ss mann and Ognyan Kounchev}
\maketitle
\begin{abstract}
In the present paper we will introduce a new approach to multivariate
interpolation by employing polyharmonic functions as interpolants, i.e. by
solutions of higher order elliptic equations. We assume that the data arise
from $C^{\infty}$ or analytic functions in the ball $B_{R}.$ We prove two main
results on the interpolation of $C^{\infty}$ or analytic functions $f$ in the
ball $B_{R}$ by polyharmonic functions $h$ of a given order of polyharmonicity $p.$


\end{abstract}

\section{Introduction and statement of results}

Interpolation theory is one of the oldest and most classical subjects of
mathematical analysis. It has been established in the work of Newton, Lagrange
and numerous other mathematicians. Interpolation plays a fundamental role in
algebraic geometry and numerical analysis, in particular in approximation of
integrals (quadrature and cubature formulas), in finite element methods, and others.

There is a number of approaches to multivariate interpolation which are based
on multivariate polynomials and radial basis functions (RBF), see e.g.
\cite{sobolev}, \cite{deboorron}, \cite{schaback}, \cite{gasca1},
\cite{olver}. From the practical point of view the problem of interpolation of
scattered data has been treated successfully by means of tools such as RBF
(see e.g. \cite{schaback}, \cite{jetter} and references there) or polysplines
(see \cite{okbook}), which are in general not globally analytic.

So far there remains the fundamental problem from the point of view of
mathematical analysis to construct a multivariate interpolation theory based
on \emph{globally analytic} tools. The multivariate polynomials fail to
deliver such tools. Indeed, it is well known and quite clear that multivariate
polynomial interpolation differs in important ways from its univariate
counterpart. The main difference is the fact that the multivariate polynomials
fail to constitute a Chebyshev system, cf. \cite{akhiezerApproximation},
\cite{karlinstudden}.

Furthermore, let us recall that in the one-dimensional case the polynomial
interpolation is closely related to a wide class of quadrature formulas. And
the existing multivariate interpolation theories mentioned above do not
provide a satisfactory theory of multivariate cubatures.

On the other hand, objects like the solutions of elliptic PDEs, in particular
the polyharmonic functions, have entered the scene of approximation and spline
theory (see e.g. \cite{kounchev85}, \cite{kounchev86}, \cite{kounchev91},
\cite{kounchev92}, \cite{haymankorenblum}, \cite{okbook} and the references
given there), and they satisfy a generalized definition of a Chebyshev system,
see \cite{kounchevChebyshevsystems}. Is there an interpolation theory based on
solutions of elliptic PDEs which provides a satisfactory analog to the
classical one-dimensional results?

In the present paper we address the above question by considering an
interpolation theory based on polyharmonic functions. Let us recall that a
function $h$ is polyharmonic of order $p$ in a domain $D\subseteq
\mathbb{R}^{n}$ if it satisfies the equation $\Delta^{p}h\left(  x\right)
=0\text{ in }D,$ cf. \cite{aron}, \cite{sobolev}. It is important to emphasize
the fact that in order to obtain satisfactory interpolation results one has to
reconsider the whole paradigm of ''set of interpolation points''. In
particular, in view of the fact that the space of polyharmonic functions is
infinite-dimensional, one may consider interpolation sets $\Gamma$ which are
the union of hypersurfaces in $\mathbb{R}^{n}.$ Some results towards this
interpolation theory have been obtained in \cite{atakhodzhaev},
\cite{kounchev85}, \cite{kounchev86}, \cite{haymankorenblum}. Let us focus on
the analogy with the one-dimensional case: one is seeking such sets $\Gamma$
which would correspond to the usual $N$ points $\left\{  x_{1},x_{2}%
,...,x_{N}\right\}  $ in $\mathbb{R}^{1}$ where a polynomial $P$ of degree
$\leq N-1$ solves the interpolation problem
\[
P\left(  x_{j}\right)  =c_{j}\qquad\text{for }j=1,2,...,N
\]
for arbitrary data $c_{j}.$ In particular, $P\left(  x_{j}\right)  =0$ for
$j=1,2,...,N$ implies $P\equiv0.$ It is clear that the main problem is to
identify multivariate analogs to the ''unisolvent'' sets $\left\{  x_{1}%
,x_{2},...,x_{N}\right\}  $.

Let us draw the reader's attention to the obstacles faced by the usual theory
of interpolation with polyharmonic functions, related to the zero sets of
polyharmonic functions. In \cite{atakhodzhaev}, \cite{kounchev86},
\cite{haymankorenblum}, and references therein, such sets of interpolation
$\Gamma\subseteq D$ have been considered which are unions of $N$ concentric
spheres. It has been proved in these works that $\Gamma$ is a set of
uniqueness, i.e. if $\Delta^{N}h=0$ in $D$ and if $h\left(  x\right)  =0$ for
all $x\in\Gamma$ then $h\equiv0$ in $D.$ So far, attempts to consider sets
$\Gamma$ with a slightly more general geometry have led to a dead-end. In
\cite{atakhodzhaev} (see the Russian edition of $1985$) Atakhodzhaev has
constructed a set of two closed convex curves $\gamma_{1}$ and $\gamma_{2}$ in
$\mathbb{R}^{2}$ with $\gamma_{1}$ contained in the convex hull of $\gamma
_{2},$ such that there exists a (non-trivial) biharmonic function $h$ with
$\Delta^{2}h=0$ inside $\gamma_{2}$ and $h\left(  x\right)  =0$ for all
$x\in\gamma_{1}\cup\gamma_{2}.$ This result has been dealt with in
\cite{haymankorenblum} as well.

The last fact completely destroys any hope of finding reasonable unisolvent
sets living in the space $\mathbb{R}^{n}.$ In the present paper we formulate a
concept of interpolation where the unisolvent sets live in what we call a
''semi-frequency domain'' which arises from the Laplace-Fourier spherical
harmonic expansion of a function, see formula (\ref{laplaceFourier}) below.

In order to motivate our approach to \emph{ polyharmonic interpolation } let
us recall that in the classical one--dimensional interpolation theory error
estimates are proved when data $c_{j}$ are obtained from a differentiable
function, i.e.
\[
c_{j}=f\left(  x_{j}\right)  \qquad\text{for }j=1,2,...,N,
\]
with $f\in C^{N+1}$. In that case one may consider estimates of the error of
interpolation
\[
E_{N}\left[  f\right]  \left(  x\right)  =f\left(  x\right)  -P_{N}\left(
x\right)  ,
\]
see \cite{krylov}, \cite{davis}. More subtle results are obtained when $f$ is
an analytic function and $N\longrightarrow\infty.$

Now let us turn to the multivariate situation. Corresponding to the univariate
case, in order to obtain a reasonable \emph{multivariate polyharmonic
interpolation theory} we will assume that the multivariate data arise from
$C^{\infty}$ or analytic functions.

Let us first introduce some necessary notions and notations. We will work in
the ball $B_{R}$ defined by
\[
B_{R}:=\left\{  x\in\mathbb{R}^{n}:\left|  x\right|  <R\right\}  .
\]
Assume that we have a basis of the space of \emph{harmonic homogeneous
polynomials} of degree $k$ (called spherical harmonics) which are denoted as
$Y_{k\ell}\left(  x\right)  $ for $k=0,1,...,$ and $\ell=1,2,...,d_{k},$
where
\begin{equation}
d_{k}=\frac{1}{\left(  n-2\right)  !}\left(  n+2k-2\right)  \left(
n+k-3\right)  \cdot\cdot\cdot\left(  k+1\right)  , \label{dk}%
\end{equation}
see \cite{steinweiss}. They are assumed to be orthonormalized with respect to
the scalar product
\[
\frac{1}{\omega_{n-1}}\int_{\mathbb{S}^{n-1}}u\left(  \theta\right)  v\left(
\theta\right)  d\theta
\]
on the unit sphere, where $\omega_{n-1}$ is the area of the unit sphere
$\mathbb{S}^{n-1}$ in $\mathbb{R}^{n}$; we have put
\[
x=r\theta,\qquad r=\left|  x\right|  .
\]

Let us denote by $C^{\infty}\left(  \overline{B_{R}}\right)  $ the set of
$C^{\infty}$ functions on a neighborhood of $\overline{B_{R}}.$ For $f\in
C^{\infty}\left(  \overline{B_{R}}\right)  $ we have the expansion in
spherical harmonics
\begin{equation}
f\left(  x\right)  =\sum_{k=0}^{\infty}\sum_{\ell=1}^{d_{k}}\widetilde
{f}_{k,\ell}\left(  r\right)  Y_{k,\ell}\left(  \theta\right)  .
\label{laplaceFourier}%
\end{equation}
We will use the following representation of $C^{\infty}$ and of analytic
functions in the ball, see \cite[p. 501, Proposition 1]{baouendi}:

\begin{proposition}
\label{PBaouendi}Let $f$ be in $C^{\infty}\left(  \overline{B_{R}}\right)  .$
Then we have the following expansion
\begin{equation}
f\left(  x\right)  =\sum_{k=0}^{\infty}\sum_{\ell=1}^{d_{k}}f_{k,\ell}\left(
r^{2}\right)  r^{k}Y_{k,\ell}\left(  \theta\right)
\label{expansionLaplaceFourier}%
\end{equation}
where the functions $f_{k,\ell}\in C^{\infty}\left(  \left[  0,R^{2}\right]
\right)  .$ The function $f$ is analytic in some neighborhood of $0$ in
$\mathbb{R}^{n}$ if and only if there exist $t_{0}>0$ and $M>0$ such that, for
all indices $k\geq0,$ $1\leq\ell\leq d_{k},$ and $j\geq0,$ we have
\begin{equation}
\sup_{0\leq t\leq t_{0}}\left|  \left(  \frac{d^{j}}{dt^{j}}\right)
f_{k,\ell}\left(  t\right)  \right|  \leq M^{k+j+1}j!.\label{EstimateBaouendi}%
\end{equation}
I.e. $f$ is analytic if and only if (\ref{EstimateBaouendi}) holds, and in
that case each function $f_{k,\ell}$ is also analytic.
\end{proposition}

On the other hand if $h$ is a function polyharmonic of order $N$ in the ball
$B_{R}$, then we have as in (\ref{expansionLaplaceFourier}) the expansion
\begin{equation}
h\left(  x\right)  =\sum_{k=0}^{\infty}\sum_{\ell=1}^{d_{k}}h_{k,\ell}\left(
r^{2}\right)  r^{k}Y_{k,\ell}\left(  \theta\right)  ,
\label{polyharmonicexpansion}%
\end{equation}
and it is well known (see Sobolev \cite{sobolev}, or \cite{okbook}) that the
coefficients $h_{k,\ell}\left(  \cdot\right)  $ are polynomials of degree
$N-1.$ Thus we may put into correspondence the functions $f_{k,\ell}$ and the
polynomials $h_{k,\ell},$ which is the core of the polyharmonic interpolation.

The polyharmonic interpolation problem is now very natural to formulate:
Assume that for every fixed $\left(  k,\ell\right)  $ with $k=0,1,...$ and
$\ell=1,2,...,d_{k}$ we have interpolation points which we assume to be
pairwise different:
\[
0\leq r_{k,\ell,1}<r_{k,\ell,2}<...<r_{k,\ell,N}\leq R.
\]
Then for every $\left(  k,\ell\right)  $ we find the polynomials $h_{k,\ell}$
of degree $\leq N-1$ from the one--dimensional interpolation problems
\begin{equation}
h_{k\ell}\left(  r_{k,\ell,j}^{2}\right)  =f_{k\ell}\left(  r_{k,\ell,j}%
^{2}\right)  \qquad\text{for }j=1,2,...,N. \label{interpolation}%
\end{equation}

Now \textit{the main question} is: For which distribution of the points
$\left\{  r_{k,\ell,j}\right\}  $ and for which functions $f$ is the series in
(\ref{polyharmonicexpansion}) convergent? If we have convergence then we will
call the function $h$ a \textbf{polyharmonic interpolant of order }$N.$ Our
first result says that for every distribution of the points $\left\{
r_{k,\ell,j}\right\}  $ and for a wide class of $C^{\infty}$ functions $f$ we
have convergence. Indeed, we have the following amazing result.

To make our result more transparent we will introduce the following seminorms
denoted by $\left\|  \cdot\right|  _{N}$, which are motivated by
(\ref{EstimateBaouendi}):%
\begin{equation}
\left\|  f\right|  _{N}:=\underset{k,\ell}{\overline{\lim}}\sup_{0\leq t\leq
R}\left|  \frac{1}{N!}\left(  \frac{d^{N}}{dt^{N}}\right)  f_{k,\ell}\left(
t\right)  \right|  ^{\frac{1}{k+N+1}}. \label{normnew}%
\end{equation}
We see that
\[
\sup_{N\geq0}\left\|  f\right|  _{N}=M
\]
where $M$ is the constant in (\ref{EstimateBaouendi}).

\begin{theorem}
\label{T1} Let the function $f$ be $C^{\infty}$ in a neighborhood of the
closed ball $\overline{B_{R}}$ and the interpolation knots $\left\{
r_{k,\ell,j}\right\}  _{k,\ell,j}$ satisfy%
\[
0\leq r_{k,\ell,j}\leq R,
\]
where $k=0,1,2,...,$ $\ell=1,2,...,d_{k}$ and $j=1,2,...,N.$ If the seminorm
$\left\|  f\right|  _{N}$ satisfies
\begin{equation}
R\left\|  f\right|  _{N}<1, \label{RM}%
\end{equation}
then there exists a unique polyharmonic interpolation function $h\left(
r\theta\right)  $ of order $N$ which belongs to $L_{2}\left(  \mathbb{S}%
^{n-1}\right)  $ for every $r\leq R,$ and $h$ belongs to $L_{2}\left(
B_{R}\right)  .$

Assuming (\ref{RM}), the error of interpolation is given by
\[
\left\|  f\left(  r\theta\right)  -h\left(  r\theta\right)  \right\|
_{L_{2}\left(  \mathbb{S}^{n-1}\right)  }\leq CR^{2N}\left\|  f\right|
_{N}^{N+1}.
\]
\end{theorem}

We see that in a certain sense the above Theorem \ref{T1} presents a complete
analog to the one--dimensional interpolation since we may take arbitrary knots
of interpolation $r_{k,\ell,j}.$ However we see that condition (\ref{RM}) is a
restriction on the arbitrariness of the data $f$ and this is the price which
we have to pay for the infinite-dimensionality of the problem. Does this
restriction imply a specialization in the one-dimensional case? The answer is
''no''. Indeed, since the one-dimensional polyharmonic functions of order $N$
are just polynomials of degree $\leq2N-1$ we see that condition (\ref{RM}) is
trivially fulfilled due to $\overline{\lim}$ in (\ref{normnew}).

There is still another way to consider the one-dimensional case embedded into
the multivariate case, namely, when in the sums (\ref{expansionLaplaceFourier}%
) and (\ref{polyharmonicexpansion}) only the term for $k=0$ is non-zero. Then
$f\left(  x\right)  =f_{0,1}\left(  r^{2}\right)  $ and $h\left(  x\right)
=h_{0,1}\left(  r^{2}\right)  $ where $h_{0,1}\left(  \cdot\right)  $ is a
polynomial of degree $\leq N-1.$ Indeed, in the one-dimensional a $C^{\infty}$
function $f$ is identical with the univariate analytic function $f_{0,1}%
\left(  \cdot\right)  $ in the expansion (\ref{laplaceFourier}), and the knots
are $r_{0,1,j}$ with $1\leq j\leq N.$ We see that in this case restriction
(\ref{RM}) is always satisfied, i.e. Theorem \ref{T1} extends the
one-dimensional theory in a natural way.

If we change the point of view, and consider $f$ to be fixed, then we have to
choose a radius $R$ small enough to fulfill (\ref{RM}).

As a second result we consider the special case of the knots which are lying
on $N$ concentric spheres in $\mathbb{R}^{n},$ i.e. when the knots $\left\{
r_{k,\ell,j}\right\}  _{k,\ell,j}$ satisfy
\[
r_{k,\ell,j}=r_{j}\qquad\text{for }j=0,1,2,...,N-1
\]
for all indices $\left(  k,\ell\right)  $. Assume that $f$ is a function
analytic in a neighborhood of $\overline{B_{R}},$ and that the polyharmonic
function $h$ is an interpolant of $f$, i.e. satisfies (\ref{interpolation}).
From the expansions in spherical harmonics (\ref{expansionLaplaceFourier}),
(\ref{polyharmonicexpansion}) for every fixed $r$, and for $j=0,1,2,...,N-1,$
we see that the interpolation problem (\ref{interpolation}) is equivalent to
the following polyharmonic interpolation problem on concentric spheres
\begin{equation}
h\left(  r_{j}\theta\right)  =f\left(  r_{j}\theta\right)  \qquad\text{for
}\theta\in\mathbb{S}^{n-1}. \label{interpolationSpheres}%
\end{equation}

Let us recall the following result from \cite[Theorem XI.3]{sobolev},

\begin{proposition}
Let $\varphi$ be a function defined and continuous on the unit sphere. A
necessary and sufficient condition for the analyticity of $\varphi$ is that in
the representation
\[
\varphi\left(  \theta\right)  =\sum_{k=0}^{\infty}\sum_{\ell=1}^{d_{k}}%
\varphi_{k,\ell}Y_{k,\ell}\left(  \theta\right)
\]
the coefficients $\varphi_{k,\ell}$ have exponential decay, i.e. there exist
two constants $K$ and $\eta>0,$ such that
\begin{equation}
\left|  \varphi_{k,\ell}\right|  \leq Ke^{-\eta k}\qquad\text{for every
}k=0,1,2,...;\ \ell=1,2,...,d_{k}. \label{analytic_Sobolev}%
\end{equation}
\end{proposition}

Let us put
\[
\varphi^{j}\left(  \theta\right)  :=f\left(  r_{j}\theta\right)
\qquad\text{for }\theta\in\mathbb{S}^{n-1}.
\]
From the estimate (\ref{analytic_Sobolev}) we see that for all
$j=0,1,2,...,N-1$ we have a number $\eta_{j}>0$ such that
\begin{equation}
\left|  \varphi_{k,\ell}^{j}\left(  \theta\right)  \right|  \leq Ke^{-\eta
_{j}k}. \label{fiklj}%
\end{equation}

Now we have again the question of convergence of the series
(\ref{polyharmonicexpansion}) and it is solved by the second main result of
our paper:

\begin{theorem}
\label{T2} Let the numbers $r_{j}$ with $0<r_{1}<r_{2}<...<r_{N}\leq R$ be
given, and for the parameters $\eta_{j}$ of the analytic functions
$\varphi^{j}$ defined in (\ref{fiklj}) the inequality
\begin{equation}
R\cdot\max_{j}\left(  \frac{e^{-\eta_{j}}}{r_{j}}\right)  <1
\label{condition_Rer}%
\end{equation}
be satisfied. Then the polyharmonic function of order of polyharmonicity $N$
satisfying the interpolation problem (\ref{interpolationSpheres}) has an
$L_{2}$--convergent series in the ball $B_{R}.$
\end{theorem}

Finally, let us remark that the polyharmonic interpolation problem
(\ref{interpolation}) may be considered as embedded in a more general scheme
of interpolation theory \cite{davis} in the following way: Let us introduce
the functionals
\[
L_{k,\ell,j}\left(  f\right)  =\frac{1}{\omega_{n-1}}\int_{\mathbb{S}^{n-1}%
}f\left(  r_{k,\ell,j}\theta\right)  Y_{k,\ell}\left(  \theta\right)
d\theta.
\]
Then the polyharmonic interpolation problem (\ref{interpolation}) may be
reformulated as the problem of finding the polyharmonic function $h$
satisfying the infinite number of equations
\[
L_{k,\ell,j}\left(  h\right)  =L_{k,\ell,j}\left(  f\right)  \qquad\text{for
all }k,\ell,j.
\]

On the other hand we have a nice demonstration of the \emph{polyharmonic
paradigm} \cite{okbook} in the present situation. As we said in the
introduction, the expectation that the knots $x_{1},x_{2},...,x_{N}$ in the
one-dimensional interpolation theory will be replaced by closed surfaces
$\gamma_{j},$ $j=1,2,...,N$ in $\mathbb{R}^{n}$ in the polyharmonic
interpolation has failed. Let us consider the sets
\[
\Gamma_{j}:=\left\{  \left(  \left(  k,\ell\right)  ,\rho_{k,\ell,j}\right)
:k=0,1,2,...;\ \ell=1,2,...,d_{k}\right\}
\]
with $\rho_{k,\ell,1}<\rho_{k,\ell,2}<...<\rho_{k,\ell,N}$ . They may be
considered as a multivariate generalization of the knots $x_{1}<x_{2}%
<...<x_{N}$ in the univariate case where $x_{j}$ is replaced by $\Gamma_{j}$.
For a better understanding of the role of the sets $\Gamma_{j}$ let us make
analogy with $\mathbb{R}^{n}$ where the boundary $\partial D$ of a star shaped
domain $D$ in $\mathbb{R}^{n}$ (centered at the origin $0$) can be written in
spherical coordinates as
\[
\partial D=\left\{  \left(  \theta,\rho_{\theta}\right)  :\text{for all
}\theta\in\mathbb{S}^{n-1}\right\}
\]
for some function $\rho_{\theta}\geq0$ defined on the sphere $\mathbb{S}%
^{n-1}.$ The results of the present paper show that the knots of interpolation
$x_{1}<x_{2}<...<x_{N}$ in the one-dimensional interpolation theory have been
replaced by the sequence of monotonely increasing ''sphere-like'' sets
$\Gamma_{1},$ $\Gamma_{2},$ ..., $\Gamma_{N}.$

\section{Proof of Theorem \ref{T1}: Polyharmonic Interpolant for General Knots}

Here we provide the proof of Theorem \ref{T1}.%

%TCIMACRO{\TeXButton{Proof}{\proof}}%
%BeginExpansion
\proof
%EndExpansion
By the definition of $h$ in (\ref{polyharmonicexpansion}) and
(\ref{interpolation}) $h_{k,\ell}$ are polynomials of degree $\leq N-1$ and we
may apply the classical results about the remainder of the interpolation,
hence
\begin{equation}
R_{N-1}\left(  x\right)  =\frac{\omega\left(  x\right)  }{N!}f^{\left(
N\right)  }\left(  \xi\right)  , \label{RN-1}%
\end{equation}
see \cite{davis} or \cite[(3.2.10)]{krylov}, and we obtain the formal series
\[
f\left(  x\right)  -h\left(  x\right)  =\sum_{k=0}^{\infty}\sum_{\ell
=1}^{d_{k}}Y_{k\ell}\left(  \theta\right)  r^{k}\frac{\omega_{k\ell}\left(
r^{2}\right)  }{N!}f_{k\ell}^{\left(  N\right)  }\left(  \xi_{k\ell}\right)
;
\]
here as usually $\omega_{k,\ell}\left(  r^{2}\right)  =\prod_{j=1}^{N}\left(
r^{2}-r_{k,\ell,j}^{2}\right)  .$

By the definition of $\left\|  f\right|  _{N}$ it follows by a standard
argument that the $L_{2}$ norm of the above is estimated by
\begin{align}
\left\|  f\left(  r\theta\right)  -h\left(  r\theta\right)  \right\|
_{L_{2}\left(  \mathbb{S}^{n-1}\right)  }^{2}  &  =\int_{\mathbb{S}^{n-1}%
}\left|  f\left(  r\theta\right)  -h\left(  r\theta\right)  \right|
^{2}d\theta\label{f-h}\\
&  \leq\sum_{k=0}^{\infty}\sum_{\ell=1}^{d_{k}}\left|  r^{k}\omega_{k\ell
}\left(  r^{2}\right)  \left\|  f\right|  _{N}^{k+N+1}\right|  ^{2}\\
&  \leq CR^{4N}\left\|  f\right|  _{N}^{2N+2}\sum_{k=0}^{\infty}\sum_{\ell
=1}^{d_{k}}\left|  r^{k}\left\|  f\right|  _{N}^{k}\right|  ^{2}.\nonumber
\end{align}

The convergence of the last series follows from the assumption (\ref{RM}).
Hence follows the $L_{2}$--convergence of the series for the polyharmonic
function $h.$ Also the estimate for the error of interpolation follows
directly.
%TCIMACRO{\TeXButton{EndProof}{\endproof}}%
%BeginExpansion
\endproof
%EndExpansion

\section{Proof of Theorem \ref{T2}: Polyharmonic Interpolation on $N$
Concentric Spheres}

Next we prove Theorem \ref{T2}.%

%TCIMACRO{\TeXButton{Proof}{\proof}}%
%BeginExpansion
\proof
%EndExpansion
Let us write the expansion of $\varphi^{j}$ in spherical harmonics
\begin{equation}
\varphi^{j}\left(  \theta\right)  =\sum_{k=0}^{\infty}\sum_{\ell=1}^{d_{k}%
}\varphi_{k,\ell}^{j}Y_{k,\ell}\left(  \theta\right)  . \label{fijseries}%
\end{equation}

Since the polyharmonic function $h$ interpolating $\varphi^{j}$ on the sphere
of radius $r_{j}$ has the form
\[
h\left(  x\right)  =\sum_{k,\ell}Y_{k\ell}\left(  \theta\right)  r^{k}%
h_{k\ell}\left(  r^{2}\right)  ,
\]
where $h_{k\ell}$ are polynomials of degree $N-1,$ we see that for all
$k=0,1,2,...,$ and $\ell=1,2,...,d_{k}$, and for all $j=0,1,2,...,N-1$ we need
to have
\[
r_{j}^{k}h_{k\ell}\left(  r_{j}^{2}\right)  =\varphi_{k\ell}^{j};
\]
hence,
\[
h_{k\ell}\left(  r_{j}^{2}\right)  =\frac{\varphi_{k\ell}^{j}}{r_{j}^{k}}.
\]
We have to prove that the series for $h$ is $L_{2}$--convergent i.e.
\[
\sum_{k,\ell}\int_{0}^{R}\left|  r^{k}h_{k\ell}\left(  r^{2}\right)  \right|
^{2}dr<\infty.
\]

First we will find estimates for all $h_{k,\ell}.$ We have the explicit
representation for the polynomials $h_{k\ell}$ in the form given in Krylov
\cite[p. 42]{krylov} and Davis \cite[p. 33]{davis}, where we put
$x_{j}=r_{k,\ell,j}^{2}.$ Let us put for the Lagrange fundamental functions
\[
\omega_{j}^{k,\ell}\left(  r^{2}\right)  :=\frac{\left(  r^{2}-x_{0}\right)
\cdot\cdot\cdot\left(  r^{2}-x_{j-1}\right)  \left(  r^{2}-x_{j+1}\right)
\cdot\cdot\cdot\left(  r^{2}-x_{N-1}\right)  }{\left(  x_{j}-x_{0}\right)
\cdot\cdot\cdot\left(  x_{j}-x_{j-1}\right)  \left(  x_{j}-x_{j+1}\right)
\cdot\cdot\cdot\left(  x_{j}-x_{N-1}\right)  }.
\]
Then we have
\begin{equation}
h_{k\ell}\left(  r^{2}\right)  =\sum_{j=0}^{N-1}\omega_{j}^{k,\ell}\left(
r^{2}\right)  \frac{\varphi_{k\ell}^{j}}{r_{j}^{k}}. \label{hinterpolation}%
\end{equation}
Bearing in mind that $r_{0}<r_{1}<\cdot\cdot\cdot<r_{N-1},$ we obtain (with
the same $K$ for all $j$'s and $\eta_{j}$'s ), the following estimate
\begin{align*}
\left|  h_{k\ell}\left(  r^{2}\right)  \right|   &  \leq K\sum_{j=0}%
^{N-1}\left|  \omega_{j}^{k,\ell}\left(  r^{2}\right)  \right|  \frac
{e^{-\eta_{j}k}}{r_{j}^{k}}\\
&  \leq\frac{K_{1}}{\delta}R^{2N}\sum_{j=0}^{N-1}\frac{e^{-\eta_{j}k}}%
{r_{j}^{k}};
\end{align*}
here $K_{1}>0$ is a suitable constant and $\delta:=\min_{j=1,2,...,N-1}\left(
\left|  x_{j}-x_{j-1}\right|  \right)  .$ Hence we obtain the estimate
\begin{align}
\sum_{k=0}^{\infty}\sum_{\ell=1}^{d_{k}}\int_{0}^{R}\left|  r^{k}h_{k\ell
}\left(  r^{2}\right)  \right|  ^{2}dr  &  \leq\sum_{k=0}^{\infty}\sum
_{\ell=1}^{d_{k}}\int_{0}^{R}r^{2k}\left|  \frac{K_{1}}{\delta}R^{2N}%
\sum_{j=0}^{N}\frac{e^{-\eta_{j}k}}{r_{j}^{k}}\right|  ^{2}dr\nonumber\\
&  \leq\left(  \frac{K_{1}}{\delta}R^{2N}\right)  ^{2}\sum_{k=0}^{\infty}%
\sum_{\ell=1}^{d_{k}}\left|  \sum_{j=0}^{N}\frac{e^{-\eta_{j}k}}{r_{j}^{k}%
}\right|  ^{2}\frac{R^{2k+1}}{2k+1}\nonumber\\
&  \leq\left(  \frac{K_{2}}{\delta}R^{2N}\right)  ^{2}\sum_{k=0}^{\infty
}k^{n-2}\left|  \sum_{j=0}^{N}\frac{e^{-\eta_{j}k}}{r_{j}^{k}}\right|
^{2}\frac{R^{2k+1}}{2k+1}. \label{inequal}%
\end{align}
To obtain the last inequality we have used the estimate
\[
d_{k}\leq Ck^{n-2}%
\]
for some constant $C>0$ which follows from (\ref{dk}). Putting
\[
M=\max_{j}\left(  \frac{e^{-\eta_{j}}}{r_{j}}\right)
\]
we obtain the estimate
\[
\sum_{j=0}^{N}\frac{e^{-\eta_{j}k}}{r_{j}^{k}}\leq C\left(  N+1\right)
M^{k}.
\]
The convergence of the series in (\ref{inequal}) follows from the assumption%
\[
R\cdot\max_{j}\left(  \frac{e^{-\eta_{j}}}{r_{j}}\right)  <1.
\]
\endproof

\begin{remark}
\textit{If}%
\[
R\max_{j}\left(  \frac{e^{-\eta_{j}}}{r_{j}}\right)  >1,
\]
\textit{ then in general one may not expect that the series representing the
polyharmonic interpolant }$h$ \textit{will be convergent.} This will be shown
by the following example.
\end{remark}

\textbf{Example.} We assume that for all $j$ we have
\begin{equation}
\frac{e^{-\eta_{j}}}{r_{j}}=C,\label{Cassumption}%
\end{equation}
so that $C=\max_{j}\left(  \frac{e^{-\eta_{j}}}{r_{j}}\right)  .$ From
(\ref{hinterpolation}) it follows that
\[
h_{k\ell}\left(  r^{2}\right)  =C^{k}\sum_{j=0}^{N-1}\omega_{j}^{k,\ell
}\left(  r^{2}\right)  ,
\]
and hence
\[
\int_{0}^{R}\left|  r^{k}h_{k\ell}\left(  r^{2}\right)  \right|  ^{2}%
dr=C^{2k}\int_{0}^{R}r^{2k}\left|  \sum_{j=0}^{N}\omega_{j}^{k,\ell}\left(
r^{2}\right)  \right|  ^{2}dr.
\]
According to the basic properties of the Lagrange coefficients (see e.g.
\cite[p. 42-43]{krylov})
\[
\sum_{j=0}^{N-1}\omega_{j}^{k,\ell}\left(  r^{2}\right)  =1,
\]
so we get
\[
\int_{0}^{R}\left|  r^{k}h_{k\ell}\left(  r^{2}\right)  \right|  ^{2}%
dr=C^{2k}\int_{0}^{R}r^{2k}dr=C^{2k}\frac{R^{2k+1}}{2k+1}.
\]
Finally, for a suitable constant $C_{1}>0$ the inequality
\[
\sum_{k,\ell}\int_{0}^{R}\left|  r^{k}h_{k\ell}\left(  r^{2}\right)  \right|
^{2}dr\geq C_{1}\sum_{k=0}^{\infty}k^{n-2}C^{2k}\frac{R^{2k+1}}{2k+1}%
\]
holds true and the last series is divergent since $CR>1.$ The proof is
finished using assumption (\ref{Cassumption}).

If assumption (\ref{Cassumption}) does not hold then we can see by standard
asymptotics arguments that for large $k$ we will have
\[
h_{k\ell}\left(  r^{2}\right)  \approx C^{k}\sum_{j=0}^{N-1}\omega_{j}%
^{k,\ell}\left(  r^{2}\right)  ,
\]
and hence
\[
\int_{0}^{R}\left|  r^{k}h_{k\ell}\left(  r^{2}\right)  \right|  ^{2}dr\geq
C_{2}C^{2k}\frac{R^{2k+1}}{2k+1}%
\]
for a suitable $C_{2}>0.$ This proves the divergence of the series.%
%TCIMACRO{\TeXButton{EndProof}{\endproof}}%
%BeginExpansion
\endproof
%EndExpansion

\textbf{Acknowledgment. } The authors acknowledge the support of the
Institutes Partnership Project V-Koop-BUL/1014793 by the Alexander von
Humboldt--Foundation, and the second author thanks the support of the
Greek--Bulgarian IST contract for the period 2005--2007.

We would like to thank the anonimous referee for several helpful remarks which
greatly enhanced the readability of the paper.

Werner Hau\ss mann, Department of Mathematics, University of Duisburg--Essen,
Lotharstr. 65, 47057 Duisburg, Germany;

haussmann@math.uni-duisburg.de

Ognyan Kounchev, Institute of Mathematics and Informatics, Bulgarian Academy
of Sciences, Acad. G. Bonchev Str. 8, 1113 Sofia, Bulgaria;

kounchev@math.bas.bg, kounchev@math.uni-duisburg.de
\end{document}